\documentclass[12pt]{amsart}
\newtheorem{lem}{Lemma}[section]
\newtheorem{lemma}{Lemma}[lem]

\newtheorem{theorem}[lem]{Theorem}
\newtheorem{corollary}[lem]{Corollary}
\theoremstyle{definition}
\newtheorem{definition}[lem]{Definition}

\newtheorem{remark}[lem]{Remark}

\newcommand{\ep}{\varepsilon}
\newcommand{\Proof}{\underbar{Proof}{\hskip 0.1in}}
\newcommand{\rn}{\mathbb{R}^N}
\newcommand{\hf}{\frac{1}{2}}
\newcommand{\qt}{\frac{1}{4}}

\newcommand{\dvol}{{\rm d}vol}
\newcommand{\dn}{{\rm d}^N}

\begin{document}
\title{PERTURBATION OF DOMAIN: SINGULAR RIEMANNIAN METRICS} 
\author{C. Mason}

\begin{abstract}
We develop further some aspects of the spectral theory of a class of Riemannian manifolds introduced by E. B. Davies; in particular we study the best constant in the Hardy Inequality which has become important in spectral theory.  One application of this constant has been to a certain type of domain perturbation.  This technique is useful when the domain has an irregular boundary and this is the case with the manifolds under consideration.  However, in this paper we show that the manifolds possess a Hardy constant that lies outside the range permitted by existing theorems.  Yet we are still able to prove theorems which give information about the domain perturbation problem and moreover, we set up a specific example which can be used to show that our results are the best possible.
\end{abstract}
\subjclass{35P99, 47A75, 47B25, 58J99}
\keywords{boundary decay, Laplacian, Hardy inequality, spectral convergence, singular Riemannian metric}
\maketitle

\section{Introduction}

In the paper {\em Two-Dimensional Riemannian Manifolds with Fractal Boundary}, \cite{RMFB}, Davies introduced a class of Riemannian manifolds which exhibited unusual properties.  These manifolds are simple in that the sense that there is only one coordinate patch and the geometry when viewed with the natural Euclidean metric is extremely straightforward (and much of the corresponding spectral theory well understood).  However, when equipped with a certain Riemannian metric, which becomes singular at the boundary, the manifolds show a link with Euclidean domains that have fractal boundaries.  More precisely Davies constructs a piecewise smooth homeomorphism between the unit disc with a specific example of the singular metric and the von Koch snowflake domain equipped with the Euclidean metric.

One strategy for computing the Dirichlet eigenvalues of a region with irregular boundary is to use the following perturbation of domain argument.  Given a Riemannian manifold $\Omega$ one considers, for small enough $\ep>0$, a smaller domain $\Omega_\ep$ such that
\[
\{ x \in \Omega : \text{dist}(x,\partial \Omega)>\ep \} \subseteq \Omega_\ep \subseteq \Omega.
\]
We denote by $H_D$ and $H_{\ep, D}$ the Dirichlet Laplacians associated with $\Omega$ and $\Omega_\ep$ respectively.  If $\lambda_{n,\ep}$ denotes the $n$th eigenvalue of $H_{\ep,D}$ and $\lambda_n$ denotes the $n$th eigenvalue of $H_D$ then variational arguments imply that
\[
\lambda_{n,\ep}\geq \lambda_n
\]
for all $n$.  Our goal is to find estimates of the quantity $\lambda_{n,\ep}-\lambda_n$ to assess how quickly these approximations converge as $\ep \to 0+$.

Results relating to this problem are given by Davies in \cite{SBEO} and by Pang in \cite{P96} and \cite{P97}.  Pang's results are applicable to the bottom eigenvalue for open, bounded and  simply connected sets.  Davies' results are valid for all eigenvalues and it is the techniques in this paper that we will follow (Pang's proofs are completely unrelated).   The key assumption made in \cite{SBEO} is that a Hardy inequality of the type
\[
c^2(H_D+a) \geq d^{-2}
\]
holds in the sense of quadratic forms.  Here $d(x):=$dist$(x,\partial \Omega)$ and $c,a$ are some constants.  The infinimum of all possible values for $c$ such that there exists an $a<\infty$ for which the above inequality holds is known as the Hardy constant.

The key result in \cite{SBEO} is the following:
\begin{align}\label{eqn:basic1}
\lambda_{n,\ep}-\lambda_n \leq c_n\ep^{\frac{2}{c}}
\end{align}
where $c_n$ is some constant and $c$ is the Hardy constant associated with the region which is assumed to satisfy
\[
c\geq 2.
\]
The restriction on $c$ is important for the proof presented in \cite{SBEO}.  In Theorem \ref{thm:riemannianhardy} we prove that the class of singular manifolds introduced in \cite{RMFB} has a Hardy constant that satisfies
\[
1<c \leq 2.
\]
This result stands in contrast to the situation of regular Euclidean domains.  It is known that the Hardy constant for every Euclidean domain whose boundary has at least one regular point is never smaller than 2.  A precise formulation of this statement and proof may be found in \cite{HC}.

Despite this we prove in Theorem \ref{thm:main} that for any manifold with Hardy constant $1\leq c$ inequality (\ref{eqn:basic1}) is still true.  The key difference is the use of better quadratic form and operator norm inequalities in Theorem \ref{thm:sbeonorm}.  The conclusions as applied to the particular Riemannian manifolds are stated explicitly in Corollary \ref{cor:main}.  We emphasise that the rate of convergence is faster than $O(\ep)$ when $c<2$.

Furthermore, in  section 7 we give a concrete example.  This is a rotationally invariant domain and after separating variables one can reduce to a one dimensional problem.  The techniques for dealing with the class of problems that arises are presented in another paper, \cite{M2}, since they are based on techniques from ordinary differential equations and are very different from anything presented here.  They can be used to show that the constant achieved in our general method is sharp for this particular example.

The proofs require two distinct notions of the distance to the boundary and the relationship between them.  The first distance is the usual Euclidean distance to the boundary and this is denoted by $\sigma$.  Since the boundary will be smooth in the Euclidean sense this is simply the length of the shortest straight line that intersects the boundary orthogonally.  This function is only guaranteed to be smooth in a neighbourhood of the boundary but in general it satisfies
\[
|\nabla_E \sigma |_E \leq 1
\]
in the distributional sense (we use the subscript here to stress that everything is the usual Euclidean quantity).

The second distance is the distance induced by the Riemannian metric.  This will be denoted $d$.  This may be defined as the infimum of the lengths of paths connecting the point to the bounday where the length is calculated with respect to the Riemannian arc length.  Theorem \ref{thm:distances} gives a simple relationship between these two and is heavily used. 

\section{The Basic Model}

We begin with $\Omega$ an open, bounded and connected subset of $\mathbb{R}^N$ that has a $C^2$ boundary and assign to it a Riemannian metric that becomes singular at the boundary in some controlled fashion.  In future this will be referred to as {\em the} Riemannian metric.  In terms of the arc length element and depending on some constant $\gamma$ to be specified later this is given by
\[
ds^2=\sigma^{-2\gamma}(x)(dx_1^2+\dots+dx_N^2)
\]
where
\begin{align}
\sigma(x)&:=\text{dist}_E(x,\partial \Omega). \notag
\end{align}

\begin{remark}  In the above and in all that follows a subscript $E$ will denote a quantity taken with respect to the Euclidean metric.  The absence of such a subscript means that the Riemannian metric is appropriate.  The fact that there are two metrics on the same domain is one of the principal challenges in the analysis.  
\end{remark}

\subsection{The Laplacian}\label{sect:laplacian}

The Riemannian volume element is given by
\[
\dvol=\sigma^{-N\gamma}(x)\dn x.
\]
The constant $\gamma$ is chosen to satisfy the relation
\[
0\leq N\gamma <1.
\]
$\Omega$ therefore has finite volume and is incomplete.  Moreover, the $L^p$ spaces are defined to be the spaces for which the following norms are finite:
\[
\Vert f \Vert_p^p:=\int_\Omega |f|^p \dvol=\int_\Omega |f|^p\sigma^{-N\gamma} \dn x.
\]
The norm on the tangent spaces derived from the Riemannian metric is such that
\[
|\nabla f|=\sigma^\gamma |\nabla_E f|_E.
\]

Since we have a metric we also have a Laplacian and this can be defined using quadratic form techniques.  The quadratic form is given by
\[
Q(f):=\int_\Omega |\nabla f|^2 \dvol =\int_\Omega |\nabla_E f|_E^2\sigma^{(2-N)\gamma}\dn x
\]
whenever this is finite. 
\begin{definition}
We define $W^{1,2}(\Omega)$ to be the space of weakly differentiable functions on $\Omega$ for which
\[
Q(f)+\Vert f \Vert_2^2<\infty
\]
and define $W^{1,2}_0(\Omega)$ to be the closure of $C^\infty_c(\Omega)$ in $W^{1,2}(\Omega)$.
\end{definition}

 \begin{remark}
If $\gamma=0$ then these are the standard Sobolev spaces which we will denote by $ W^{1,2}_E(\Omega)$ and $ W^{1,2}_{0,E}(\Omega)$
\end{remark}

\begin{definition}
The Dirichlet Laplacian $H_D$ is associated in the standard way with the quadratic form $Q$ with domain $W^{1,2}_0(\Omega)$.  Similarly the Neumann Laplacian $H_N$ is found from $Q$ with the domain $W^{1,2}(\Omega)$.

The usual Dirichlet and Neumann Laplacians (whose forms are defined on $W^{1,2}_{0,E}(\Omega)  $ and $W^{1,2}_E(\Omega)  $ respectively) will be denoted by $-\Delta_D$ and $-\Delta_N$.

\end{definition} 

Formally, the differential operators are given by
\[
Hf:=-\sigma^{N\gamma} \sum_i \frac{\partial}{\partial x_i}\left( \sigma^{(2-N)\gamma}\frac{\partial f}{\partial x_i}\right)
\]
and in particular for $N=2$
\[
Hf=-\sigma^{2\gamma}\Delta_E f
\]
although these are for illustration only and the rigorous definitions above are to be preferred.

\subsection{Manifolds with Fractal Boundary}

The final section of \cite{RMFB} gives the construction of a piecewise smooth homeomorphism between the ball equipped with a metric of the above type and $S$, the Koch snowflake domain,  equipped with the Euclidean metric such that the two metrics are Lipschitz equivalent (i.e. quasi-isometric).  It is commented in this section that the technique seems to be restricted to two dimensions.  This result allows one to establish  a bound for the heat kernel $K(t,x,y)$ of $\Omega$ with the Riemannian metric,  subject to Neumann conditions:
\[
0<K(t,x,y) \leq k t^{-1} \text{ for } 0<t<1
\]
where $k$ is some constant.  This is possible because there is a body of knowledge about the snowflake domain.  More specifically it has the $W^{1,2}(S)\rightarrow W^{1,2}(\mathbb{R}^2)$ extension property being a quasidisc (see e.g. \cite{MS}).
\subsection{Other Examples}\label{sect:11}

In \cite{P88} Pang considers singular elliptic operators in a similar situation.  The operator acts in $L^2(\Omega, \eta^\lambda(x)\dn x)$ and is derived in the usual way from the quadratic form

\[
Q_1(f)=\int_\Omega \eta^\mu|\nabla_Ef|_E^2\dn x \text{ where  } f\in C_c^\infty(\Omega).
\]
$\Omega$ is again a smooth domain in $\rn$, $\lambda$ and $\mu$ are constants and $\eta$ is a smooth function such that $\eta(x) \sim \sigma(x)$ as $x \rightarrow \partial \Omega$.  The proofs are based on differential operator techniques  and he systematically establishes a variety of $L^1$ and $L^2$ spectral properties such as conservation of probability and essential self-adjointness.  The key point is that these properties depend only on the values of $\lambda$ and $\mu$.  It should also be mentioned that as well as the above case he also considered the case $\Omega=\rn$ and $\eta \sim (1+|x|^2)^{\frac{1}{2}}$ and in this case he was later able to weaken the smoothness assumptions on $\eta$ in \cite{P95}.  For results related to this case see also \cite{L}.

\section{Riemannian Distance to the Boundary}\label{sect:regdist}

We denote the Riemannian distance between two points $x,y\in\Omega$ by $d(x,y)$.  The Riemannian distance can be calculated in two distinct ways which we recall below.  The first allows us to compute upper bounds.

\begin{definition}
A smooth curve $w:[\alpha,\beta]\rightarrow\Omega$ is said to be regular if $w'(t)\ne0$ for $t\in (\alpha,\beta)$.  A continuous map $w:[\alpha,\beta]\rightarrow\Omega$ is called piecewise regular if there exists a finite subdivision $\alpha=a_0<a_1<\cdots<a_k=\beta$ such that $w|_{[a_{i-1},a_{i}]}$ is regular for $i=1,\dots,k$. Piecewise regular curves will be known as {\em admissible curves}
\end{definition}
\begin{definition}
If the Riemannian arc length is given by
\[
ds^2:= g_{ij}(x)dx_idx_j
\]
then an admissible curve $w:=[\alpha,\beta]\rightarrow\Omega$ with coordinates $w^{i}(t)$ has length $\ell (w)$ defined by
\[
\ell(w):=\int_\alpha^\beta \sqrt{ g_{ij}\frac{dw^i}{dt}\frac{dw^j}{dt}  }dt
\]
using the summation convention.  Now we define $d(x,y)$ by
\[
d(x,y):=\inf \{ \ell(w) : w \text{ is an admissible curve and }w(\alpha)=x \text{ and }w(\beta)=y\}.
\] 
\end{definition}
It can be shown that
\[
d(x,y)=\sup \left\{ |\psi(x) -\psi(y)|:\psi \text{ is Lipschitz continuous and } |\nabla \psi | \leq 1 \right\}.
\]
This can give useful lower bounds.

Also considered in \cite{RMFB} is the completion of $\Omega$ and it gives formulae for the geodesics connecting two points $v_1,v_2\in \partial \Omega$.  These are curves in $\Omega$ that intersect $\partial \Omega$ at right angles.  One further important piece of notation is the following:
\[
d(x):=\text{dist}(x,\partial \Omega)
\]
for $x\in\Omega$. As usual
\[
\text{dist}(x,\partial \Omega):=\inf_{y\in\partial\Omega}d(x,y).
\]
An important result is the relationship between $d$ and $\sigma$.  

\begin{theorem}\label{thm:distances}
Let $\Omega$ and $\gamma$ be as in the basic model.  For $x\in \Omega$ we have
\[
d(x)=\frac{\sigma^{1-\gamma}(x)}{1-\gamma}
\]
\end{theorem}

\Proof
Let $\psi$ be given by
\[
\psi_0(x) = \frac{\sigma^{1-\gamma}(x)}{1-\gamma}
\]

A calculation then shows that in the sense of distributions

\begin{align}
|\nabla_E \psi_0|_E^2&=
\sigma^{-2\gamma}|\nabla_E \sigma |_E^2 \notag \\
&\leq \sigma^{-2\gamma} \notag
\end{align}
and
\[
|\nabla \psi_0|^2=\sigma^{2\gamma} |\nabla_E \psi_0 |_E^2 \leq 1.
\]

Now, $\psi_0\in W^{1,\infty}$ and so $\psi_0$ is Lipschitz continuous, satisfies $|\nabla \psi_0| \leq 1$  and hence
\[
d(x)  \geq \frac{\sigma^{1-\gamma}(x)}{1-\gamma}.
\]

Next, given a point $x\in\Omega$, let $\Lambda$ be the straight line segment that minimises the Euclidean distance to the boundary, i.e. the straight line segment that orthogonally intersects the boundary at the nearest point in the Euclidean sense.  Choose the coordinate system so that the $x_1$ axis lies along $\Lambda$ and the origin is at the point of intersection of $\Lambda$ and $\partial \Omega$.  In this coordinate system we have $x=(a,0,\dots,0)$ for some $a$.

Now for any point $y=(t,0,\dots,0)$ lying along the straight line segment $\Lambda$ we have  $\sigma(y)=t$.  Thus, the length of $\Lambda$ is given by
\[
\int_0^{a}t^{-\gamma}dt=\frac{{a}^{1-\gamma}}{1-\gamma}= \frac{\sigma^{1-\gamma}(x)}{1-\gamma}.
\]

Hence
\[
d(x) \leq \frac{\sigma^{1-\gamma}(x)}{1-\gamma}
\]
and we have the result by combining the two inequalities.
\section{The Hardy Inequality}\label{sect:hi}

\subsection{Generalised Hardy Inequalities}

The Hardy Inequality is a useful tool when investigating the spectral properties of elliptic operators.  It will play a role in what follows but first we stop to gather the essential definitions and results.

We begin with $M$ a Riemannian manifold, $L$ a non-negative, second order elliptic operator acting in $L^2(M, \dvol)$ subject to Dirichlet boundary conditions and $\rho$ a positive, continuous function defined on $M$ which satisfies $|\nabla \rho| \leq 1$ at least in the distributional sense ($\rho$ will be taken to be the distance to the boundary).  Also let $Q(f)$ be the quadratic form associated with $L$.

\begin{definition}
We say that $L$ satisfies a {\em weak Hardy inequality} with respect to $\rho$ if there exists a constant $c>0$ and a constant $a <\infty$ such that 
\[
\int_\Omega \frac{|f|^2}{\rho^2} \dvol \leq c^2 \left( Q(f) + a \Vert f \Vert^2 \right)
\]
for all $f$ in the domain of $Q$.  The weak Hardy constant is then the infimum of all possible $c$ such that this holds. 

If we may take $a=0$ in the above then we say $L$ satisfies a {\em strong Hardy Inequality}, and the infimum is then called the strong Hardy constant.  
\end{definition}

Further information may be found in the review by Davies (\cite{ROHI}) - also of note is the book by Opic and Kufner (\cite{OK}) which establishes conditions for a Hardy Inequality to hold in a region with H\"older class boundary.  No attempt is made to quantify the constant in general.  

We conclude these introductory remarks with an important standard result.  For a proof see \cite{ROHI}.

\begin{theorem}\label{thm:classhardy}
If $\Omega \subset \mathbb{R}^N$ (usual Euclidean metric) is bounded with a $C^2$ boundary and $\Delta_D$ is the associated Dirichlet Laplacian then there exists $a\in\mathbb{R}$ such that
\[
\sigma^{-2} \leq 4(-\Delta_D+a)
\]
in the sense of quadratic forms.  
\end{theorem}
\subsubsection{HI for the Basic Model}

We need the following result which is trivial to prove but an important observation.

\begin{lemma}\label{lem:basic}
\[
\int_\Omega \frac{|f(x)|^2}{d(x)^2}\dvol=(1-\gamma)^2\int_\Omega \frac{|f(x)|^2}{\sigma(x)^2}\sigma^{(2-N)\gamma}(x)\dn x
\]
if either side is finite.
\end{lemma}

\Proof Simply recall that $\dvol=\sigma^{-N\gamma}\dn x$ and $d(x)=\frac{\sigma^{1-\gamma}(x)}{1-\gamma}$.

Thus if we could establish a Hardy Inequality of the form
\begin{align}\label{eqn:basic}
\int_\Omega \frac{|f|^2}{\sigma^2}\sigma^{(2-N)\gamma}\dn x \leq c^2 \left( \int_\Omega |\nabla_E f|_E^2 \sigma^{(2-N)\gamma}\dn x+a\int_\Omega |f|^2 \sigma^{(2-N)\gamma}\dn x \right)
\end{align}
then we would immediately have
\[
\int_\Omega \frac{|f|^2}{d^2}\dvol \leq c^2(1-\gamma)^2\left( \int_\Omega |\nabla_E f|_E^2 \sigma^{(2-N)\gamma}\dn x+a'\int_\Omega |f|^2 \sigma^{-N\gamma}\dn x \right)
\]
where $a'=a\| \sigma^{2\gamma} \|_\infty$ and $f\in W^{1,2}_0$.  

To prove inequality (\ref{eqn:basic}) we work with the Hilbert space $L^2(\Omega,\sigma^{(2-N)\gamma}\dn x)$ and use the following:

\begin{theorem}\cite[Theorem 1.5.12]{HKST}\label{thm:opweak}
Suppose $L$ is an elliptic operator on $L^2(\Omega)$ and that there is a positive continuous function $\phi\in W^{1,2}_{\text{loc}}(\Omega)$ and a potential $V$ in $L^1(\Omega)$ such that
\[
L\phi \geq V\phi
\]
in the distributional sense.  Then we have
\[
L \geq V
\]
in the sense of quadratic forms.
\end{theorem}

Now we can prove the first main theorem.

\begin{theorem}\label{thm:hardy}
If $\Omega$ is convex (in the Euclidean sense) and equipped with the Riemannian metric then the strong Hardy Inequality holds with strong constant $c$ satisfying
\[
c \leq \frac{2(1-\gamma)}{1+(N-2)\gamma}.
\]
\end{theorem}

\Proof  Let $L$ denote the Friedrichs extension of the  operator given formally by
\[
Lf:=-\sigma^{(N-2)\gamma}\sum_i  \frac{\partial}{\partial x_i} \left(\sigma^{(2-N)\gamma} \frac{\partial f}{\partial x_i} \right).
\]
$L$ is defined in terms of quadratic forms as usual.  Now, let $\phi$ be defined by
\[
\phi:=\sigma^\alpha
\]
where
\[
\alpha:=\frac{1+(N-2)\gamma}{2}.
\]
Then we have
\begin{align}
L\phi&=-\alpha(\alpha-1+(2-N)\gamma)\sigma^{\alpha-2} |\nabla_E \sigma |_E^2- \alpha \sigma^{\alpha-1}\Delta_E \sigma  \notag \\
&\geq -\alpha(\alpha-1+(2-N)\gamma)\sigma^{\alpha-2} \notag \\
&=\frac{(1+(N-2)\gamma)^2 }{4}\sigma^{-2}\phi \notag
\end{align}
in the sense of distributions.
This calculation uses the fact that $\Delta_E \sigma \leq 0$ for any convex set $\Omega$.

Now, applying Lemma \ref{thm:opweak} we have
\[
\int_\Omega \frac{|f|^2}{\sigma^2}\sigma^{(2-N)\gamma}\dn x \leq \frac{4}{(1-(2-N)\gamma)^2} \int_\Omega |\nabla_E f|^2\sigma^{(2-N)\gamma}\dn x
\]
and applying Lemma \ref{lem:basic} gives
\[
\int_\Omega \frac{|f(x)|^2}{d(x)^2}\dvol\leq \frac{4(1-\gamma)^2}{(1-(2-N)\gamma)^2} \int_\Omega |\nabla f|^2\dvol.
\]

\begin{theorem}\label{thm:riemannianhardy}
The weak Hardy Inequality holds for all $\Omega$ in the basic model with constant $c$ satisfying
\[
c \leq \frac{2(1-\gamma)}{1+(N-2)\gamma}.
\]
\end{theorem}

\Proof  Let $\Phi$ be a smooth function that equals $\sigma$ near to the boundary.  We now consider the differential operator $K$ defined on $C^2_c$ functions say, by
\[
Kf:=-\Phi^{(N-2)\gamma} \frac{\partial}{\partial x_i} \left(\Phi^{(2-N)\gamma} \frac{\partial f}{\partial x_i} \right)
\]
and again we denote by $K$ its Friedrichs extension.

The trial function $\phi$ is now chosen to be a smooth function which near the boundary satisfies
\[
\phi=\sigma^\alpha-\sigma
\]
where $\alpha$ is as before.

Then,
\begin{align}
\frac{K\phi}{\phi}&=\frac{-\alpha(\alpha-1+(2-N)\gamma)\sigma^{-2}-\alpha\sigma^{-1}\Delta_E\sigma+(2-N)\gamma\sigma^{-1-\alpha}+\sigma^{-\alpha}\Delta_E\sigma}{1-\sigma^{1-\alpha}} \notag \\
&=\frac{(1+(N-2)\gamma)^2}{4}\sigma^{-2}+(\frac{(1+(N-2)\gamma)^2}{4}-(N-2)\gamma)\sigma^{-1-\alpha}+O(\sigma^{-2\alpha}) \notag
\end{align}
provided we are close enough to the boundary.

Now observe that 
\[
\frac{(1+(N-2)\gamma)^2}{4}-(N-2)\gamma=\frac{(1-(N-2)\gamma)^2}{4}>0
\]
to conclude that
\[
\frac{K\phi}{\phi} \geq \frac{(1+(N-2)\gamma)^2}{4}\sigma^{-2}
\]
near the boundary.  The behaviour of $\phi$ away from the boundary is not important for the existence of the weak Hardy Inequality and we conclude that
\[
K\phi(x) \geq \frac{(1+(N-2)\gamma)^2}{4}\sigma^{-2}(x)\phi(x)-a\phi(x)
\]
for all $x\in\Omega$ provided $a$ is chosen to be large enough.  An application of Theorem \ref{thm:opweak} establishes the final result.

\begin{remark}
The previous theorem indicates the general point that the weak Hardy Inequality depends only on the local geometry of the boundary.  To gain information about the strong version we must make assumptions about the global geometry. 
\end{remark}

There is a link between the strong Hardy constant and the Minkowski dimension of the boundary.  This link allows us to prove that the result given above is sharp.

The following definition is taken from \cite{ROHI} - it is not the usual definition but is sufficient for our purposes.

\begin{definition} We say that the boundary $\partial \Omega$ has interior Minkowski dimension $\alpha$ if there exist positive constants $c_1$ and $c_2$ such that
\[
c_1 \ep^{N-\alpha} \leq \text{vol}(\partial \Omega_\ep) \leq c_2 \ep^{N-\alpha}
\]
where
\[
\partial \Omega_\ep:=\{ x \in \Omega : \text{dist}(x,\partial \Omega)<\ep \} 
\]

for small enough $\ep>0$ (vol$(A)$ denotes the Riemannian volume of the set $A$.  Note also that the distance is calculated with respect to the Riemannian metric).
\end{definition}

\begin{theorem}\label{min}\cite[Theorem 3.1]{RMFB}
Let $\alpha$ denote the interior Minkowski dimension of $\partial \Omega$.  Then
\[
N-1\leq \alpha = \frac{N-1}{1-\gamma} < N
\]
\end{theorem}

The following result of Davies and Mandouvalos is now of interest:

\begin{theorem} \label{min2} \cite[Theorem 6]{ROHI}
If $\partial \Omega$ has interior Minkowski dimension $\alpha$, so that $\alpha \geq N-1$, then c, the strong Hardy constant of $\Omega$ with respect to the Laplacian, satisfies
\[
c(2+\alpha - N) \geq 2
\]
\end{theorem}

In conclusion:
\begin{theorem}
Let $\Omega$ be convex and let $c$ denote the strong Hardy constant.  Then
\[
c=\frac{2(1-\gamma)}{1+(N-2)\gamma}.
\]
\end{theorem}

\Proof Combine Theorems \ref{thm:hardy}, \ref{min} and \ref{min2}.

The range of values that $c$ may take are of importance and we make a note of this now.

\begin{lemma}
For $\gamma \in [0,1/N)$ the number \[c(\gamma)=2(1-\gamma)/(1+(N-2)\gamma)\] satisfies
\[
1<c\leq 2
\]
\end{lemma}

\Proof Recall that 
\[
c(\gamma)=\frac{2(1-\gamma)}{1+(N-2)\gamma}
\]
and $0\leq \gamma <1/N$.  

A straightforward calculation shows that $c'(\gamma)<0$ (where $'$ indicates differentiation with respect to $\gamma$), $c(0)=2$ and $c(1/N)=2(N-1)/(2N-2)=1 $.

\section{Sobolev Spaces and Norm Estimates}

There are two Hilbert Spaces which have featured thus far: the usual Euclidean space $L^2(\Omega,\dn x)$ which has norm
\[
\Vert f \Vert_{2,E}^2=\int_\Omega |f(x)|^2\dn x
\]
and the space $L^2(\Omega,\dvol)$ which is derived from the Riemannian metric and has the norm
\[
\Vert f \Vert_2^2=\int_\Omega |f(x)|^2\dvol=\int_\Omega |f(x)|^2\sigma^{-N\gamma}\dn x.\]

\begin{lemma}
$L^2(\Omega,\dvol)\subset L^2(\Omega,\dn x)$
Moreover we have
\[
\Vert \sigma^{\hf N\gamma}\Vert_\infty^{-1}\Vert f \Vert_{2,E} \leq \Vert f \Vert_2.
\]
\end{lemma}
\Proof
This follows from the fact that
\[
\Vert \sigma^{\hf N\gamma} f \Vert_{2}=\Vert f \Vert_{2,E}.
\]

Turning to the quadratic forms we see that the usual Euclidean quadratic form
\[
Q_E(f)=\int_\Omega |\nabla_E f|_E^2\dn x
\]
and the Riemannian form 
\[
Q(f)=\int_\Omega |\nabla f|^2\dvol =\int_\Omega |\nabla_E f|_E^2\sigma^{(2-N)\gamma}\dn x
\]
are equal when $N=2$.  Moreover, we can see from the previous lemma that $Q_E(f)\leq Q(f)$ in general.  We now turn to the Sobolev spaces defined previously (see section \ref{sect:laplacian}).  We recall first a definition.

\begin{definition}
We say $\Omega$ is a regular domain if a generalised Hardy Inequality holds: i.e. there exist constants $c_1,c_2$ with $c_1>0$ such that
\[
H_D\geq c_1/d^2-c_2
\]
where $H_D$ is the natural Dirichlet Laplacian and $d$ the distance to the boundary.
\end{definition}
\begin{remark}Theorems \ref{thm:classhardy} and \ref{thm:riemannianhardy} imply that $\Omega$ is a regular domain with both the Euclidean and the Riemannian metric.
\end{remark}
The following theorem gives a complete description of the Sobolev space $W^{1,2}_0$.

\begin{theorem}\label{thm:reg}
If $\Omega$ is a regular domain then
\[
W^{1,2}_0(\Omega)=W^{1,2}(\Omega)\cap \left\{f:\int_\Omega \frac{|f|^2}{d^2}<\infty \right\}.
\]
\end{theorem}

\Proof  This theorem is essentially proved in \cite[Theorem 1.5.6]{HKST}.  The proof involves a regularised distance function.  

\begin{lemma}
If $f\in W^{1,2}_{0,E}(\Omega)$ then $f\in L^2(\Omega,\dvol)$.  Moreover,
\[
\Vert f \Vert_2=\Vert \sigma^{-\hf N\gamma}f\Vert_{2,E} \leq \Vert f \sigma^{-1}\Vert_{2,E} \Vert \sigma^{1-\hf N\gamma}\Vert_\infty.
\]
\end{lemma}

\Proof   If $f\in W^{1,2}_{0,E}(\Omega)$ then Theorem \ref{thm:reg} ensures that
\[
\int_\Omega \frac{|f|^2}{\sigma^2} {\rm d}^Nx <\infty
\]
i.e. $\sigma^{-1}f\in L^2(\Omega,\dn x)$.  Also $\sigma^{1-\hf N\gamma}\in L^\infty$ and so H\"{o}lder's Inequality immediately gives us that 
\[
\Vert \sigma^{-N\gamma/2}f\Vert_{2,E} \leq \Vert f \sigma^{-1}\Vert_{2,E} \Vert \sigma^{1-\hf N\gamma}\Vert_\infty .
\]

\begin{lemma}\label{sobolev}
For arbitrary dimension $N$ we have
\[
W^{1,2}(\Omega) \subset  W^{1,2}_E(\Omega)
\]
and hence $W^{1,2}_0(\Omega) \subset  W^{1,2}_{0,E}(\Omega)$.

Moreover, if $N=2$ we have
\[
W^{1,2}_0(\Omega)= W^{1,2}_{0,E}(\Omega).
\]
\end{lemma}

\Proof The first two statements follow immediately from the discussion about the norms and the quadratic forms.    Now, if $f\in W^{1,2}_{0,E}(\Omega)$ then  $f\in L^2(\Omega,\dvol)$ by previous lemma.  Next observe $Q$ is independent of $\gamma$ when $N=2$ and so it follows that $f\in W^{1,2}(\Omega)$.  Finally we observe that\[
\int_\Omega \frac{|f|^2}{d^2}\dvol = (1-\gamma)^2   \int_\Omega \frac{|f|^2}{\sigma^2} {\rm d}^2x <\infty
\]
and thus we conclude that $f\in W^{1,2}_0(\Omega)$.

\section{Boundary Perturbation}

We now return to the boundary perturbation problem as described in the introduction.  We recall that in the case that a weak Hardy Inequality holds in the domain with Hardy constant $c\geq 2$ the result is as follows:
\[
\lambda_{n,\ep}-\lambda_n = O(\ep^{2/c}).
\]

The difficulty in applying these results to our situation is that $1<c\leq 2$ and indeed $c<2$ whenever $\gamma>0$.  Of course, the Hardy Inequality remains true with the larger constant $c=2$ and we could simply apply these results to get the fact that
\[
\lambda_{n,\ep}-\lambda_n = O(\ep).
\]
This is less than satisfactory.  Numerical experiments and specific examples suggest a rate of convergence which is faster than this - in the following section we set up an example that may be used to show a faster rate.

The results in this section will be general in character although they will ultimately be applied to the basic model.  For this reason we think of $\Omega$ as being {\em some} Riemannian manifold with Dirichlet Laplacian $H_D$ and Hardy Inequality \[ c^2(H_D+a) \geq d^{-2}\] where $d$ is the distance to the boundary and $c>1$ some constant.

We will modify the approach in \cite{SBEO}.  This uses the original eigenvectors $\phi_i$ of $H_D$ to approximate the eigenvectors of $H_{\ep,D}$ in the variational formula.  We define a rapidly decreasing function $\mu:\Omega\rightarrow [0,\infty)$ by
\[
\mu(x):=\left\{ \begin{array}{cc} 0 & \text{ if } 0<d(x)\leq\ep \\
\frac{d(x)-\ep}{\ep} & \text{ if } \ep < d(x) \leq 2 \ep \\
1 & \text{ otherwise.} \end{array} \right.
\]
We note that $0 \leq \mu \leq 1$, $|\nabla \mu| \leq \frac{1}{\ep}$ and $\mu$ has support in $\Omega_\ep$.

The effect of applying this cut-off function is captured in the following theorem.

\begin{theorem}\label{thm:sbeogen}
For $f\in$Dom($H_D)$ the following holds
\begin{enumerate}
\item $\mu f \in $Dom($Q$),
\item $Q(\mu f) \leq Q(f) + 2\int_{\{x:\ep<d(x)<2\ep\} }|\nabla f|^2 \dvol+2\ep^{-2}\int_{\{x:\ep<d(x)<2\ep\} }|f|^2\dvol$,
\item $\| f\|_2 \geq \| \mu f\| \geq \| f\|_2-\sqrt{\int_{\{x:0<d(x)<2\ep\} }|f|^2\dvol}$.
\end{enumerate}
\end{theorem}

\Proof Theorem \ref{thm:sbeogen}:1  More generally we have the following :\begin{lemma}\cite[Lemma 2]{SBEO}\label{lem:wdom}
If $f\in$Dom($Q$) and $w\in W^{1,\infty}$ then $wf\in$Dom($Q$)
\end{lemma}

\Proof[Lemma \ref{lem:wdom}] \begin{align}
Q(wf)&=\int_\Omega (|w\nabla f+f\nabla w|^2)\dvol \notag \\
&\leq \int_\Omega (2w^2|\nabla f|^2+2|f|^2|\nabla w|^2)\dvol \notag \\
&\leq 2\Vert w \Vert_\infty^2Q(f)+2\Vert \nabla w\Vert_\infty^2\Vert f\Vert_2^2. \notag 
\end{align}

Now we may give the conclusion of Theorem \ref{thm:sbeogen}:1.

It has already been noted that $\mu\in W^{1,\infty}$ and so the lemma can be applied.

\Proof Theorem \ref{thm:sbeogen}:2    Let $S:=\{x : \ep < d(x) < 2 \ep\}$. Then we have
\begin{align}
Q(\mu f) -Q(f) \leq& \int_S |\nabla (\mu f)|^2 \dvol\notag \\
\leq& 2\int_S \mu^2|\nabla f|^2\dvol+2\int_S |\nabla \mu|^2|f|^2\dvol \notag \\
\leq& 2\int_S |\nabla f|^2\dvol+2\ep^{-2}\int_S|f|^2\dvol. \notag
\end{align}

\Proof Theorem \ref{thm:sbeogen}:3\begin{align}
| \|f \|_2- \|\mu f \|_2|^2 &\leq \| f-\mu f\|_2^2 \notag \\
&=\int_\Omega (1-\mu)^2|f|^2 \dvol \notag \\
&\leq \int_{\{x:0<d(x)<2\ep\}} |f|^2\dvol . \notag
\end{align}
This completes the proof.

The key point now is to obtain boundary decay estimates by considering the following integral
\[
\int_{\{x:d(x)<\ep\}}|f|^2 \dvol
\]
making only the assumption that a Hardy Inequality of the form
\[
d^{-2}\leq c^2(H_D+a)
\]
holds in the sense of quadratic forms for $c>1$.   

Although we have relaxed the condition on $c$ this is not without incurring some further cost.  The estimates achieved in \cite{SBEO} were valid for functions lying in the domain of the operator.  The following theorem captures the key results of \cite{SBEO}.

\begin{theorem}\label{thm:goal}
Let $f \in $Dom($H_D$), $c^2(H_D+a) \geq d^{-2}$ for  $c\geq 2$.  Then
\[
\int_{\{x:d(x)<\ep\}}|f|^2\dvol \leq c^{2+2/c}\| (H_D+a)f\|_2 \|(H_D+a)^{1/c}f\|_2
\]
and
\[
\int_{\{x:d(x)<\ep\}}|\nabla f|^2\dvol \leq (c^{2/c}+c^{2/c}(1+c)^{2+2/c})\| (H_D+a)f\|_2 \|(H_D+a)^{1/c}f\|_2
\]
\end{theorem}

\Proof See \cite[Theorems 4 and 8]{SBEO}

Our results will be even more restrictive in that they demand that the functions lie in the space Dom$((H_D+a)^{1/2+1/c})$.  This is not a serious restriction when the results are applied to the perturbation of domain problem.  In this we are only concerned with the eigenvectors of $H_D$ which lie in Dom($(H_D+a)^s)$ for all $s \geq 0$.  Indeed, when the results are extended to more general elliptic operators with measurable coefficients the distinction becomes perhaps even less important because in this case it is extremely difficult to determine even Dom$(H_D)$ and certainly it cannot usually be identified with one of the standard Sobolev spaces.

The following result is taken from \cite{SBEO}, it does not depend  critically on the value of $c$ and therefore is also applicable in our situation.

\begin{definition}
We define the function $w$ by
\[ w(x):=(\max\{d(x),\ep\})^{-1/c}.
\]
\end{definition}

\begin{lemma}\cite[Lemma 3]{SBEO}\label{thm:sbeoold}
If $f\in$Dom($H_D$) then
\[
\int_\Omega \left(\frac{w^2}{c^2 d^2}-|\nabla w|^2\right) |f|^2\dvol \leq \langle Hf,w^2f\rangle+a\Vert wf\Vert_2^2.
\]
\end{lemma}

\Proof By the Hardy Inequality and Lemma \ref{lem:wdom} we have $wf\in$Dom($Q$) and
\[
\int_\Omega \frac{w^2f^2}{c^2 d^2}\dvol \leq Q(wf)+a\Vert wf\Vert_2^2.
\]
Moreover,
\[
Q(wf)=\langle H_Df,w^2 f\rangle+\int_\Omega |\nabla w|^2|f|^2\dvol.
\]
This calculation may be found in slightly more detail in \cite{SBEO} but is routine.

The main theorem in this section differs from a similar result in \cite{SBEO} in that we can ignore the restriction that $c\geq 2$.

\begin{theorem}\label{thm:sbeonew}
Let $f\in$ Dom$((H_D+a)^{1/2+1/c})$ where $1\leq c \leq 2$. Then
\[
|\langle H_Df,w^2f\rangle +a\Vert wf\Vert_2^2| \leq c^{2/c}\|(H+a)^{1/2+1/c}f\|_2\|(H_D+a)^{1/2}f\|_2
\]
\end{theorem}

\Proof
We have
\[
w^{2c} \leq d^{-2} \leq c^2(H_D+a)
\]
in the sense of quadratic forms.  Thus
\[
(H_D+a)^{-\hf}w^{2c}(H_D+a)^{-\hf} \leq c^2
\]
and so
\[
\| w^c (H_D+a)^{-\hf} \| \leq c.
\]
We also have
\begin{align}\label{eqn:quadest}
w^{2(2-c)} \leq (d^{-2})^{\frac{2-c}{c}} \leq (c^2(H_D+a))^{\frac{2-c}{c}}
\end{align}
since $0 \leq (2-c)/c \leq 1$.  This uses the fact that for non-negative self adjoint operators $A$ and $B$ the quadratic form inequality $0 \leq A \leq B$ implies $0 \leq A^\alpha \leq B^\alpha$ for all $0<\alpha<1$.  It is at this point that the assumption $c\geq 1$ is important.  For a proof of this fact see \cite[Lemma 4.20]{OPS}.

Now from (\ref{eqn:quadest}) we also have the fact that
\[
\| w^{2-c}(H_D+a)^{-\frac{2-c}{2c}} \|\leq c^{\frac{2-c}{c}}.
\]
Finally we make the definition $H_a:=H_D+a$ and then
\begin{align*}
|\langle H_Df,w^2 f \rangle +a\Vert wf\Vert_2^2|&=|\langle (H+a)f,w^2f\rangle|  \\
&=|\langle w^{2-c}H_a f, w^c f\rangle| \\
&=|\langle w^{2-c}H_a^{-\frac{2-c}{2c}}H_a^{\frac{2-c}{2c}}H_af,w^cH_a^{-\hf}H_a^{\hf}f\rangle| \\
&\leq c^{\frac{2-c}{c}+1}\|H_a^{\frac{2-c}{2c}+1}f\|_2\|H_a^{\hf}f\|_2 \\
&= c^{2/c}\|(H+a)^{1/2+1/c}f\|_2\|(H_D+a)^{1/2}f\|_2
\end{align*}

The theorem can now be applied in the same manner as \cite{SBEO}.  The following theorem is Theorem 4 from \cite{SBEO} using the above estimate in the proof.

\begin{theorem}\label{thm:sbeonorm}
Let $f\in$ Dom$((H_D+a)^{1/2+1/c})$. Then
\[
\int_{ \{x:d(x)<\ep \}} |f|^2 \dvol \leq c^{2+2/c}\ep^{2+2/c}\Vert (H_D+a)^{1/2+1/c} f\Vert_2\Vert (H_D+a)^{1/2}f\Vert_2
\]
\end{theorem}

\Proof Combining Theorems \ref{thm:sbeoold} and \ref{thm:sbeonew} we have
\[
\int_\Omega \left(\frac{w^2}{c^2 d^2}-|\nabla w|^2\right) |f|^2\dvol \leq  c^{2/c}\Vert (H_D+a)^{1/2+1/c} f\Vert_2\Vert (H_D+a)^{1/2}f\Vert_2
\]
If $d(x)\geq \ep$ then
\[
|\nabla w|^2 \leq c^{-2}d^{-2-2/c}=\frac{w^2}{c^2d^2}.
\]
Thus
\[
\frac{w^2}{c^2 d^2}-|\nabla w|^2 \geq 0.
\]
Also, if $d(x)<\ep$ then
\[
\frac{w^2}{c^2 d^2}-|\nabla w|^2=\frac{1}{c^2\ep^{2/c}d^2}.
\]
Hence
\begin{align}
\int_{\{x:d(x)<\ep\}} \frac{|f|^2}{d^2} \dvol&\leq c^2\ep^{2/c}\int_\Omega\left(\frac{w^2}{c^2 d^2}-|\nabla w|^2\right) |f|^2 \dvol \notag \\
&\leq c^{2+2/c}\ep^{2/c} \Vert (H_D+a)^{1/2+1/c} f\Vert_2\Vert (H_D+a)^{1/2}f\Vert_2 \notag
\end{align}

The result follows immediately.

The previous results allow us to obtain similar estimates for $\nabla f$. 
\begin{theorem}\label{thm:sbeograd}
Let $f\in$ Dom$((H_D+a)^{1/2+1/c})$. Then
\[
\int_{\{x:d(x)<\ep \}}|\nabla f|^2 \dvol \leq c_1 \ep^{2/c} \Vert (H_D+a)^{1/2+1/c} f\Vert_2\Vert (H_D+a)^{1/2}f\Vert_2.
\]
where $c_1:=c^{2/c}+c^{2/c}(1+c)^{2+2/c}$.
\end{theorem}

\Proof This is similar to Theorem 8 in \cite{SBEO}.  Again we use Theorem \ref{thm:sbeonew} in place of Theorem 1 of \cite{SBEO} at the appropriate point.

As previously mentioned, when applying the results achieved in the last section we consider only the eigenfunctions of the original operator $H_D$ and use them to approximate the eigenfunctions of the restricted operator $H_{D,\ep}$.  

\begin{definition}
Let $\phi_1,\phi_2\dots \phi_n$ denote the eigenfunctions of $H_D$, where $\langle \phi_i,\phi_j\rangle =\delta_{i,j}$, corresponding to the eigenvalues  $\lambda_1, \lambda_2, \dots, \lambda_n$.

Now, define the $n$-dimensional subspace $L_n$ by
\[
L_n:=\text{span}\{\phi_1,\dots, \phi_n\}
\]
\end{definition}

We now state a corollary to Theorems \ref{thm:sbeonorm} and \ref{thm:sbeograd}.
\begin{corollary}\label{cor:eigest}
Let $f\in L_n$.  Then
\[
\int_{\{x:d(x)<\ep\}} |f|^2\dvol \leq c^{2+2/c}\ep^{2+2/c}(\lambda_n+a)^{1+1/c}\|f\|_2^2
\]
and
\[
\int_{\{x:d(x)<\ep\}} |\nabla f|^2 \dvol \leq c_1 \ep^{2/c}(\lambda_n+a)^{1+1/c}\|f\|_2^2.
\]
\end{corollary}

\Proof It is immediate from the spectral theorem that \[\Vert (H_D+a)^pf \Vert_2 \leq (\lambda_n+a)^p\Vert f\Vert_2\] for $f\in L_n$. Now apply Theorems \ref{thm:sbeonorm} and \ref{thm:sbeograd}.

\begin{theorem}\label{thm:main}
Let $c \geq 1$.  Then, there exist constants $c_n$ such that
\[
\lambda_{n,\ep}-\lambda_n \leq c_n \ep^{2/c}
\]
for $\ep>0$ small enough.
\end{theorem}

\Proof  We use the variational estimate
\[
\lambda_{n,\ep} \leq \sup_{f \in L_n} \frac{Q(\mu f)}{\|\mu f \|_2^2}.
\]
Combining Theorem \ref{thm:sbeogen} and Corollary \ref{cor:eigest} we find that
\[
\lambda_{n,\ep} \leq \frac{\lambda_n+c'\ep^{2/c}(\lambda_n+a)^{1+1/c}} {1-c^{1+1/c}\ep^{1+1/c}(\lambda_n+a)^{\hf(1+1/c)  }}
\]
where $c'=2(c_1+c^{2+2/c})$.

\begin{corollary}\label{cor:main}
Let $\Omega$ be as given in the basic model, satisfying a Hardy Inequality with constant
\[
c=\frac{2(1-\gamma)}{1+(N-2)\gamma}.
\]
Then there exist constants $c_n$ such that 
\[
\lambda_{n,\ep}-\lambda_n \leq c_n \ep^{(1+(N-2)\gamma)/(1-\gamma)}
\]
for  $\ep>0$ small enough.
In particular if $N=2$ we have the following
\[
\lambda_{n,\ep}-\lambda_n =O(\ep^{\frac{1}{1-\gamma}}) \text{ as } \ep \to 0+.
\]
\end{corollary}

\section{The Rotationally Invariant Case}

We now consider a specific example which is particularly tractable and may be used to show that the power achieved  in the previous general result is sharp.  We use the unit disc as the underlying region, i.e.
\[
\Omega:=\{ x : |x| <1 \} \subset \mathbb{R}^2
\]
and so the Riemannian metric is given by
\[\begin{array}{cc}
ds^2=(1-|x|)^{-2\gamma}|dx|^2 & 0\leq \gamma <\frac{1}{2} \end{array}.
\]
Using the orthogonal group we make the decomposition
\[
L^2(\Omega)=\Sigma_{n=-\infty}^{\infty}L^2_n(\Omega)
\]
where
\[
L^2_n(\Omega)=\left\{f(r,\theta)=\frac{1}{\sqrt{2\pi}}g(r)e^{i n\theta}\right\}.
\]
Each subspace $L_n^2$ is invariant with respect to the Laplacian and so we can restrict our attention to each subspace.  This uses the well known process of separation of variables.  Thus we can consider only 
\begin{align}
\Vert f \Vert_2^2= \int_0^1|g(r)|^2(1-r)^{-2\gamma}rdr \notag \\
Q_n(f)= \int_0^1\left( |g'(r)|^2+\frac{n^2}{r^2}|g(r)|^2\right) rdr. \notag
\end{align}

We reiterate that the case $N=2$ causes an important simplification: the quadratic form and its domain are independent of $\gamma$.  Formally, the associated differential operator is given by

\begin{align}\label{ode:0}
Lg:=-(1-r)^{2\gamma}\left(\frac{1}{r}\frac{d}{dr}\left(r\frac{dg}{dr}\right)-\frac{n^2}{r^2}g \right) \text{ on } L^2((0,1), r(1-r)^{-2\gamma}dr).
\end{align}

Next we reformulate the problem is terms of a Schr\"odinger operator.

\begin{lemma}
Under suitable transformations we have the following expressions for the norm and quadratic form:
\begin{align}
\Vert h \Vert^2_2 &:=\int_0^\alpha |h(t)|^2dt \notag \\
Q(h)&=\int_0^\alpha \left( |h'(t)|^2+V(t)h(t) \right)dt \notag
\end{align}
where
\[
\alpha:=\frac{1}{1-\gamma}
\]
and
\[
V(t)=\frac{2\gamma-\gamma^2}{4(1-\gamma)^2}\frac{1}{(\alpha-t)^2}+\left(n^2-\qt\right)\frac{(1-\frac{t}{\alpha})^{2\gamma \alpha}}{(1-(1-\frac{t}{\alpha})^{\alpha})^2}.
\]
\end{lemma}
\Proof
Let $x=w(t)=1-((1-\gamma)(\alpha-t))^{\alpha}$ so that $t=0$ corresponds to $x=0$ and $t=\alpha$ corresponds to $x=1$.  We work with $g\in C^\infty_c(\Omega)$ and then prove the final result by approximation.  The norm and form become
\begin{align}
\Vert g \Vert_2^2&=\int_0^\alpha |g(w(t))|a(t)dt \notag \\
Q(f)&=\int_0^\alpha \left( |g'(w(t))|^2dt+|g(w(t))|b(t)\right) a(t)dt \notag
\end{align}

where
\begin{align}
a(t)&:=\frac{w(t)}{w'(t)} \notag \\
b(t)&:= \left(n\frac{w'(t)}{w(t)}\right)^2. \notag
\end{align}

Make the substitution $h(t):=g(w(t))\sqrt{a(t)}$, then

\begin{align}
 h'&=g' a^\hf +\hf g \frac{a'}{a^\hf} \notag \\
&=g' a^\hf + \hf h\frac{a'}{a}. \notag 
\end{align}

Thus
\begin{align}
\int_0^\alpha |g'|^2adt  =& \int_0^\alpha \left(|h'|^2+  \qt(\frac{a'}{a})^2|h|^2 \right)dt - \int_0^\alpha hh'\frac{a'}{a}dt \notag \\
=& \int_0^\alpha \left(|h'|^2+  \qt(\frac{a'}{a})^2|h|^2 \right)dt + \int_0^\alpha \hf |h|^2 \left( \frac{a'}{a}\right)'dt \notag \\
=& \int_0^\alpha |h'|^2 dt +\qt \int  (\frac{a'}{a})^2|h|^2dt + \hf \int_0^\alpha  |h|^2 \left( \frac{a''a-(a')^2}{a^2}\right)dt \notag\\
=&  \int_0^\alpha \left( |h'|^2 +(\frac{a''}{2a} -\qt (\frac{a'}{a})^2)|h|^2 \right) dt. \notag
\end{align}

Immediately we have
\[
V(t)=\frac{a''(t)}{2a(t)}-\frac{1}{4}\left(\frac{a'(t)}{a(t)} \right)^2+b(t).
\]

A routine, but messy, calculation establishes that in terms of $w$ we have
\[
V(t)= \frac{3}{4} \left( \frac{w''}{w'} \right)^2-\frac{w'''}{2w'}+\left( n^2-\qt \right)\left(\frac{w'}{w}\right)^2.
\]

The final form of the potential is found after utilising the formulae:

\begin{align}
w(t)=&1-\left(1-\frac{t}{\alpha}\right)^{\alpha} \notag \\
w'(t)=& \left(1-\frac{t}{\alpha}\right)^{\alpha-1} \notag \\
w''(t)=& -\gamma  \left(1-\frac{t}{\alpha}\right)^{\alpha-2} \notag \\
w'''(t)=& \gamma(2\gamma-1)\left(1-\frac{t}{\alpha}\right)^{\alpha-3}. \notag
\end{align}

After some calculations we see that we can reduce the problem to considering
\begin{align}\label{eqn:final}
-h''(t)+V(t)h(t)=\lambda h(t)
\end{align}
for
\[ 0<t<\alpha
\]
and
\[
0<t<\alpha-\ep 
\]
where
\begin{align*}
V(t) = \frac{n^2-1/4}{t^2}+O(t^{-1}) \text{ as } t \rightarrow 0+ \\
V(t) = \frac{2\gamma-\gamma^2}{4(1-\gamma)^2} \frac{1}{(\alpha-t)^2}+O((\alpha-t)^{\frac{2\gamma}{1-\gamma}})  \text{ as } t \rightarrow \alpha-.
\end{align*}

In another paper, \cite{M2}, we deal with general differential equations of this type.  The techniques and proofs are different from anything in this paper.  We refer to \cite{M2} for details and here simply quote the final result as it applies to our model.

\begin{theorem}
Let $\lambda_n$ be the nth eigenvalue of equation (\ref{eqn:final}) on the full interval $(0,\alpha)$ and $\lambda_{n,\ep}$ the $n$th eigenvalue on $(0,\alpha-\ep)$.  Then
\[
\lambda_{n,\ep}=\lambda_n+c_n\ep^{\frac{1}{1-\gamma}}+o(\ep^{\frac{1}{1-\gamma}}) \text{ as }\ep\to 0^+
\]
for some constant $c_n$.
\end{theorem}

The power in this expansion agrees with that in Corollary \ref{cor:main} and shows that it gives the best result possible.

\textbf{Acknowledgements}  I would like to thank Brian Davies for suggesting this problem and for his guidance and advice during this work.  I also acknowledge the support of the Engineering and Physical Sciences Research Council through a research studentship.

Department of Mathematics \\
King's College London \\
Strand \\
London \\
WC2R 2LS \\
U.K. \\
cmason@mth.kcl.ac.uk

\end{document}